\documentstyle{article}

\def\HDS{\vrule width0pt height2.3ex depth1.05ex\displaystyle}
\def\f#1#2{{{\HDS #1}\over{\HDS #2}}}

\def\mn{\!-\!}

\def\prop#1#2{\vspace{2ex} \noindent\textsc{#1.} \emph{#2} \par \vspace{2ex}}
\def\oprop#1#2{\noindent\textsc{#1.} \emph{#2} \par \vspace{2ex}}
\def\dkz{\noindent\textsc{Proof. }}
\def\qed{\hfill $\dashv$\vspace{2ex}}

\def\L{{\cal L}}

\begin{document}

\title{\textbf{Graphs of Plural Cuts}}
\author{\textsc{Kosta Do\v sen} and \textsc{Zoran Petri\' c}
\\[1ex]
{\small Mathematical Institute, SANU}\\[-.5ex]
{\small Knez Mihailova 36, p.f.\ 367, 11001 Belgrade,
Serbia}\\[-.5ex]
{\small email: \{kosta, zpetric\}@mi.sanu.ac.rs}}
\date{}
\maketitle

\begin{abstract}
\noindent Plural (or multiple-conclusion) cuts are inferences made
by applying a structural rule introduced by Gentzen for his
sequent formulation of classical logic. As singular
(single-conclusion) cuts yield trees, which underlie ordinary
natural deduction derivations, so plural cuts yield graphs of a
more complicated kind, related to trees, which this paper defines.
Besides the inductive definition of these oriented graphs, which
is based on sequent systems, a non-inductive, graph-theoretical,
combinatorial, definition is given, and to reach that other
definition is the main goal of the paper. As trees underlie
multicategories, so the graphs of plural cuts underlie
polycategories. The graphs of plural cuts are interesting in
particular when the plural cuts are appropriate for sequent
systems without the structural rule of permutation, and the main
body of the paper deals with that matter. It gives a combinatorial
characterization of the planarity of the graphs involved.
\end{abstract}

\vspace{.5ex}


\section{Introduction}
Plural cut is a structural inference rule introduced by Gentzen in
\cite{G35} for his plural sequent system of classical logic.  A
plural sequent (more often called multiple-conclusion sequent, or
something like that) is a sequent $\Gamma\vdash\Delta$ where
$\Delta$, as $\Gamma$, may be a collection (sequence, multiset or
set) of formulae with more than one member (see \cite{S74},
Theorem 1.2, \cite{SS78}, Chapters 1, 2, 5, \cite{G81}, Chapter
1.1, Theorem 13, and \cite{D99} for results about the relationship
between singular and plural consequence relations). Plural cut as
formulated by Gentzen with sequents based on sequences of formulae
$\Gamma$, $\Theta$, $\Delta$ and $\Lambda$ is the following rule:
\[
\f{\Gamma\vdash\Theta,A\hspace{3em} A,\Delta\vdash\Lambda}
{\Gamma,\Delta\vdash\Theta,\Lambda}
\]
A sequent $\Gamma\vdash\Delta$ is singular when the collection of formulae  $\Delta$ cannot have more than one member, and singular cut is obtained from Gentzen's plural cut by assuming that $\Theta$ is empty and that $\Lambda$ cannot have more than one member.

Gentzen assumed his rule of plural cut together with the structural rule of permutation, on both the left and right of the turnstile $\vdash$, so that the exact place of the formula $A$ in his formulation of plural cut is not essential.
Besides the plural cut rule as stated by Gentzen, the following plural cut rules:
\[
\f{\Gamma\vdash A,\Theta\hspace{3em} \Delta,A\vdash\Lambda}
{\Delta,\Gamma\vdash\Lambda,\Theta}
\]
\[
\f{\Gamma\vdash A\hspace{3em} \Delta_1,A,\Delta_2\vdash\Lambda}
{\Delta_1,\Gamma,\Delta_2\vdash\Lambda}\hspace{5em}
\f{\Gamma\vdash\Theta_1,A,\Theta_2\hspace{3em} A\vdash\Lambda}
{\Gamma\vdash\Theta_1,\Lambda,\Theta_2}
\]
were considered in \cite{A91} and \cite{L93} as appropriate for
plural sequent systems where one does not assume the structural
rule of permutation. A detailed study of cut elimination in the
context of these rules may be found in \cite{HS95}. (Something
analogous in a different area may be found in the definition of
the literal shuffle of \cite{B87}, Part~1, p.\ 29.) Let us call
these four kinds of plural cuts \emph{planar plural cuts} (as the
literature suggests).

Planar plural cuts are found in the polycategories of \cite{CS97},
which were called planar polycategories in \cite{K05}. These
polycategories differ from the polycategories of \cite{S75} where
we have the following plural cut rule:
\[
\mbox{\rm (PC)}\hspace{1em}\f{\Gamma\vdash\Theta_1,A,\Theta_2\hspace{3em} \Delta_1,A,\Delta_2\vdash\Lambda}
{\Delta_1,\Gamma,\Delta_2\vdash\Theta_1,\Lambda,\Theta_2}\hspace{3em}
\]

This rule involves a kind of permutation, which is manifested in the crossings of the following diagram:
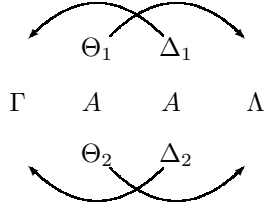
\begin{center}
\begin{picture}(90,80)

\qbezier(5,15)(30,-10)(55,15) \qbezier(35,15)(60,-10)(85,15)
\put(9,11){\vector(-1,1){5}} \put(81,11){\vector(1,1){5}}

\qbezier(5,65)(30,90)(55,65) \qbezier(35,65)(60,90)(85,65)
\put(9,69){\vector(-1,-1){5}} \put(81,69){\vector(1,-1){5}}

\put(0,40){\makebox(0,0){$\Gamma$}}

\put(30,20){\makebox(0,0){$\Theta_2$}}
\put(28,40){\makebox(0,0){$A$}}
\put(30,60){\makebox(0,0){$\Theta_1$}}

\put(60,20){\makebox(0,0){$\Delta_2$}}
\put(58,40){\makebox(0,0){$A$}}
\put(60,60){\makebox(0,0){$\Delta_1$}}

\put(90,40){\makebox(0,0){$\Lambda$}}

\end{picture}
\end{center}
These crossings require that we have the structural rule of
permutation on the left and on the right in order to state the
equations implicit in the definition of polycategory of \cite{S75}
(see P3 in Section~2; in the first of these equations, which are
analogous to the equations that stand behind our Propositions 2.2
and 2.3, we must permute $\Gamma_2$ with $\Delta_2$ and $\Gamma_3$
with $\Delta_3$, and in the second we must permute $\Delta_1$ with
$\Phi_1$ and $\Delta_2$ with $\Phi_2$). Planar plural cuts are
obtained from (PC) by requiring that either $\Theta_1$ or
$\Delta_1$ be empty and that either $\Theta_2$ or $\Delta_2$ be
empty, so that the crossings do not arise.

In this paper our main goal is to characterize in a
graph-theoretical, combinatorial, manner the planarity involved in
planar plural cuts. To achieve that, we define in three different
manners a kind of oriented graph, which we call K-graph. (The
notion of oriented graph, and other notions we need concerning
these graphs, and directed graphs in general, are defined in
Section~2.) The name of K-graphs is derived from the form of these
graphs, that may resemble up to a point a rotated K:
\begin{center}
\begin{picture}(20,8)
\put(0,0){\line(1,0){20}}

\put(10,0){\line(-3,2){10}}

\put(10,0){\line(3,2){10}}

\end{picture}
\end{center}
(see the picture below; this form resembles equally a rotated X).

Our first definition, given in Section~2, is inductive. With it,
K-graphs are obtained from some basic K-graphs by applying
operations that correspond to planar plural cuts. This definition
yields our notion of global K-graph, which is closest to planar
polycategories.

It corresponds actually to a notion somewhat more general than the
notion of planar polycategory, which we could call \emph{compass}
polycategory. Compass polycategories would be defined like planar
polycategories, but instead of having polyarrows with sources and
targets made of \emph{sequences} of objects, in compass
polycategories we would have these sources and targets made of
multisets of objects with two distinguished objects, if the
multiset is not a singleton. We refer to these distinguished
objects by $N$ and $S$ (which stand for \emph{north} and
\emph{south} respectively; we take inspiration from the compass
because in $\Gamma\vdash\Delta$ we have $\Gamma$ on the west and
$\Delta$ on the east.) In sequences of objects, the $N$ and $S$
object are the first and last object. The collections of objects
in the polyarrows of compass polycategories need not however be
sequences. We need $N$ and $S$ to characterize the operations on
global K-graphs that correspond to planar plural cuts, and we do
not need anything else. The assumption that we have sequences is
not necessary to characterize these operations.

The polyarrows of a freely generated compass polycategory may be
identified with global K-graphs where the inner vertices (a vertex
of a directed graph is inner when an edge ends in it and another
one begins in it; see Section~2) are labelled by the free
generators of the polycategory, and the remaining vertices are
labelled by objects of the polycategory. Our Propositions 2.4 and
2.5 (see Section~2) contain the essence of a completeness proof of
our notion of global K-graph with respect to compass
polycategories, which as planar polycategories are characterized
by three equations that stand behind our Propositions 2.1-2.3, and
by additional equations involving the identity polyarrow. To
simplify the exposition, we deal separately in Section~6 with
matters involving this identity. This section brings a
mathematically not very essential addition to the preceding
exposition in the main body of the paper.

It is not our intention in this paper to deal with compass
polycategories. We leave this topic for another place.

Our third definition of K-graph, given in Section~5, is non-inductive and it does not mention $N$ and $S$ any more.
It is purely graph-theoretical, and by showing the equivalence of the notion the third definition gives with the
notion of global K-graph we have achieved the main goal of the paper.

Our second definition, in Section~3, gives the notion of local K-graph, which is intermediary between the notion of global K-graph given by the first definition and the notion of K-graph given by the third definition. Our main definition of local K-graph is non-inductive as the third definition, but it still involves $N$ and $S$, as the first definition. We give however also an inductive definition of local K-graph. The notion of local K-graph given by the second definition, which is equivalent to the notion of global K-graph, as proved in Section~4, helps us to prove in Section~5 the equivalence mentioned in the preceding paragraph.

With the third definition of K-graph, the planarity of planar plural cuts, or rather their compass character, is characterized in a way that can be compared to Kuratowski's way of characterizing the planarity of graphs (see \cite{H69}, Chapter 11, and the first part of the proof of Proposition 5.5). The two approaches may be compared, but the results involved are different. In our case, we do not deal in fact with planarity, but with a related notion involving $N$ and $S$. We deal also with a special kind of oriented graph, whereas Kuratowski was concerned with the planarity of ordinary, non-directed, graphs.

The third definition yields the following picture. An arbitrary
K-graph may very roughly be described as having in the middle a
non-circular line of edges with changing directions, which we call
the transversal. Together with the transversal we have two sets of
trees, each set of a different kind: the first set consists of
trees oriented towards the root, and the second of trees oriented
towards the leafs. Both kinds of trees are planted with their
roots in the transversal. Here is an example:
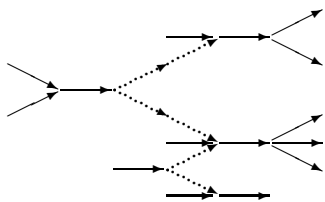
\begin{center}
\begin{picture}(130,90)
\put(0,40){\vector(2,1){19}} \put(0,60){\vector(2,-1){19}}
\put(20,50){\vector(1,0){19}} \put(40,20){\vector(1,0){19}}
\put(60,10){\vector(1,0){17.5}} \put(60,30){\vector(1,0){17.5}}
\put(60,70){\vector(1,0){17.5}} \put(80,10){\vector(1,0){19}}
\put(80,30){\vector(1,0){19}} \put(80,70){\vector(1,0){19}}
\put(100,30){\vector(1,0){19}} \put(100,29){\vector(2,-1){19}}
\put(100,31){\vector(2,1){19}} \put(100,69){\vector(2,-1){19}}
\put(100,71){\vector(2,1){19}}

\multiput(41,50)(2,1){18}{\makebox(0,0){\circle*{.5}}}
\put(74.5,67){\vector(2,1){5}} \put(55.5,57.5){\vector(2,1){5}}

\multiput(41,50)(2,-1){18}{\makebox(0,0){\circle*{.5}}}
\put(74.5,33){\vector(2,-1){5}} \put(55.5,42.5){\vector(2,-1){5}}

\multiput(61,20)(2,1){9}{\makebox(0,0){\circle*{.5}}}
\put(74.5,27){\vector(2,1){5}}

\multiput(61,20)(2,-1){9}{\makebox(0,0){\circle*{.5}}}
\put(74.5,13){\vector(2,-1){5}}

\end{picture}
\end{center}
In the middle, drawn with dotted lines, is the transversal, on the
left of which, growing westward, we have trees oriented towards
the root, and on the right of which, growing eastward, we have
trees oriented towards the leafs (for details see Section~5). The
combinatorial essence of the planarity of K-graphs is that the
transversal is non-circular (more precisely, asemicyclic; see
Section~2) and linear (more precisely, non-bifurcating; see
Section~5).

With singular cuts we would obtain just trees, oriented towards the root. One bases on such trees derivations in ordinary natural deduction, and also the notion of multicategory. A limit case of singular cut is ordinary composition in categories, which yields as graphs just chains. With K-graphs we do not have trees, but we have not gone very far away from the notion of tree.

With plural cuts in general, which are based on the cut rule (PC), we are further removed from trees, and we obtain a notion of oriented graph, which we call Q-graph, simpler to define than our notion of K-graph, both inductively and non-inductively. We investigate this notion, which when defined non-inductively reduces essentially to a weak form of connectedness and a weak form of non-circularity, in Section~7, the last section of the paper. An arbitrary Q-graph may be pictured as an arbitrary K-graph, with the transversal and two sets of trees, but the transversal is not linear any more.

The notion of global K-graph is the notion that should be used to prove by induction that every K-graph can be geometrically realized in the plane in the following special manner. A point that realizes a vertex $a$ that is not inner has the first coordinate $0$ if an edge begins in $a$, and it has the first coordinate $1$ if an edge ends in $a$. We require moreover in this realization that for every edge $(a,b)$ of our K-graph the first coordinate of the point that realizes the vertex $a$ is strictly smaller than the first coordinate of the point that realizes the vertex~$b$.

Conversely, for an oriented graph of a special kind, which is connected and non-circular in a weak sense, and satisfies moreover a condition concerning its vertices that are not inner (see conditions (1)-(3) in Section~3), we should be able to prove that if it is realized in the plane in the special manner above, then it is a K-graph. The proof of that would be inductive too, and would rely on the notion of global K-graph. We will not go here into this rather geometrical matter, which however would not improve significantly our mathematical perception of the geometrical planarity of K-graphs. We suppose that the notion of global K-graph suffices for that. The accent in this paper is put on other matters, like our third definition of K-graph, which characterizes the planarity of these graphs in a combinatorial way.

\section{Global K-graphs}
In this section we deal with our first definition of K-graph, which yields the notion of global K-graph. We establish for this notion a completeness result in Propositions 2.1-2.5, which will help us for the equivalence proofs in later sections. We start first with some elementary notions of graph theory.

A \emph{digraph} $D$ is an irreflexive binary relation on a finite nonempty set, called the set of \emph{vertices} of $D$. The ordered pairs in $D$ are its \emph{edges}. An edge $(a,b)$ \emph{begins} in $a$ and \emph{ends} in $b$.

An \emph{oriented graph} is an antisymmetric digraph.

A vertex of a digraph $D$ is a \emph{W-vertex} ($W$ stands for \emph{west}) of $D$ when in $D$ there are no edges ending in this vertex. It is an \emph{E-vertex} ($E$ stands for \emph{east}) of $D$ when in $D$ there are no edges beginning in this vertex (which means that it is a $W$-vertex of the digraph converse to $D$). It is an \emph{inner vertex} of $D$ when it is neither a $W$-vertex nor an $E$-vertex of $D$.

An edge of $D$ is a \emph{W-edge} of $D$ when it begins in a $W$-vertex of $D$, and it is an \emph{E-edge} of $D$ when it ends in an $E$-vertex of $D$. It is an \emph{inner edge} of $D$ when it begins in an inner vertex of $D$ and ends in an inner vertex of $D$.

Intuitively, in logical terms, the $W$-vertices should be
understood as labelled by premises, i.e.\ formulae from the
left-hand side of sequents, while the $E$-vertices are labelled by
conclusions, i.e.\ formulae from the right-hand side of sequents.
This is because we write from west to east. Otherwise, we could as
well understand everything in the opposite way. The inner vertices
should be understood in logical terms as corresponding to rules of
inference, i.e.\ sequents.

Throughout the paper we use $X$ as a variable standing for $W$ or $E$, and sometimes instead of $X$ we also use $Z$ for the same purpose. We assume that $\bar{W}$ is $E$ and $\bar{E}$ is $W$. We reserve the variable $Y$ for $N$ or $S$ (which stand for \emph{north} and \emph{south} respectively).

A $W$-edge $(a,b)$ of $D$ is \emph{functional} when $(a,c)\in D$ implies $b=c$. An $E$-edge $(b,a)$ of $D$ is \emph{functional} when $(c,a)\in D$ implies $b=c$ (i.e., it is functional as a $W$-edge of the digraph converse to $D$).

A \emph{basic} K\emph{-graph} $B$ is an oriented graph of the form
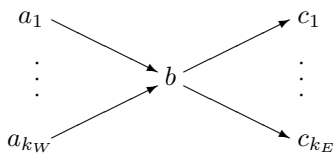
\begin{center}
\begin{picture}(110,65)
\put(10,10){\vector(2,1){40}} \put(10,55){\vector(2,-1){40}}
\put(60,35){\vector(2,1){40}} \put(60,30){\vector(2,-1){40}}

\multiput(5,26.5)(0,6){3}{\makebox(0,0){\circle*{.5}}}
\multiput(105,26.5)(0,6){3}{\makebox(0,0){\circle*{.5}}}

\put(7,55){\makebox(0,0)[r]{$a_1$}}
\put(11,8){\makebox(0,0)[r]{$a_{k_W}$}}

\put(55,37){\makebox(0,0)[t]{$b$}}

\put(103,55){\makebox(0,0)[l]{$c_1$}}
\put(103,8){\makebox(0,0)[l]{$c_{k_E}$}}

\end{picture}
\end{center}
for $k_W,k_E\geq 1$, together with the distinguished $W$-edges $NW(B)$ and $SW(B)$ and the distinguished $E$-edges $NE(B)$ and $SE(B)$, which satisfy the following condition for every $X\in\{W,E\}$:
\begin{tabbing}
\hspace{1.5em}(XYB)\hspace{3em}if $k_X\geq 2$, then $NX(B)\neq SX(B)$.
\end{tabbing}

Let $D_X$ be an oriented graph with a functional $\bar{X}$-edge $e_X$. Here $X$ can be $W$ and $E$, and we assume that $D_W$ and $D_E$ are disjoint digraphs, by which we mean that their sets of vertices are disjoint. We assume also that $e_W$ is $(a,b)$, $e_E$ is $(c,d)$ and $e$ is $(a,d)$.

Then let $D_W[e_W\mn e_E]D_E$ be the oriented graph
\[(D_W-\{e_W\})\cup(D_E-\{e_E\})\cup\{e\}\]
on the union of the vertices of $D_W$ and $D_E$ with the vertices
$b$ and $c$ omitted. This is illustrated by the following picture:
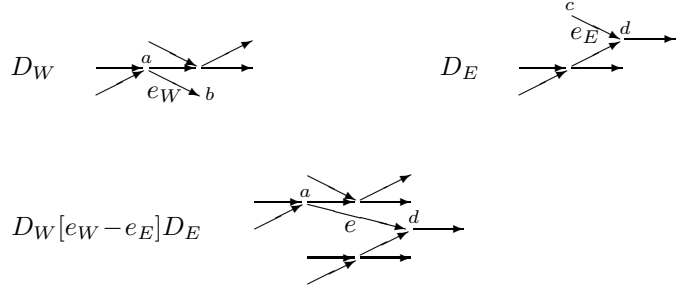
\begin{center}
\begin{picture}(200,45)
\put(0,10){\vector(2,1){18}} \put(0,20){\vector(1,0){17.5}}
\put(20,19){\vector(2,-1){19}} \put(20,20){\vector(1,0){17.5}}
\put(40,20){\vector(1,0){19}} \put(20,30){\vector(2,-1){18}}
\put(40,21){\vector(2,1){19}}

\put(-15,20){\makebox(0,0)[r]{$D_W$}}
\put(19,25){\scriptsize\makebox(0,0)[t]{$a$}}
\put(43,13){\scriptsize\makebox(0,0)[t]{$b$}}
\put(26,13){\makebox(0,0)[t]{$e_W$}}

\put(160,10){\vector(2,1){18}} \put(160,20){\vector(1,0){17.5}}
\put(180,21){\vector(2,1){18}} \put(180,20){\vector(1,0){19}}
\put(180,40){\vector(2,-1){18}} \put(200,31){\vector(1,0){19}}

\put(145,20){\makebox(0,0)[r]{$D_E$}}
\put(179,42){\scriptsize\makebox(0,0)[b]{$c$}}
\put(200,33){\scriptsize\makebox(0,0)[b]{$d$}}
\put(185,31){\makebox(0,0)[b]{$e_E$}}
\end{picture}
\end{center}
\begin{center}
\begin{picture}(80,60)
\put(0,30){\vector(2,1){18}} \put(0,40){\vector(1,0){17.5}}
\put(20,39){\vector(4,-1){37}} \put(20,40){\vector(1,0){17.5}}
\put(40,40){\vector(1,0){19}} \put(20,50){\vector(2,-1){18}}
\put(40,41){\vector(2,1){19}}

\put(20,9){\vector(2,1){18}} \put(20,19){\vector(1,0){17.5}}
\put(40,20){\vector(2,1){18}} \put(40,19){\vector(1,0){19}}
\put(60,30){\vector(1,0){19}}

\put(19,45){\scriptsize\makebox(0,0)[t]{$a$}}

\put(-20,30){\makebox(0,0)[r]{$D_W[e_W\mn e_E]D_E$}}

\put(60,32){\scriptsize\makebox(0,0)[b]{$d$}}
\put(36,33){\makebox(0,0)[t]{$e$}}

\end{picture}
\end{center}
The oriented graphs $D_W$ and $D_E$ may be conceived as obtained
from $D_W[e_W\mn e_E]D_E$ by \emph{cutting} the edge $e$ into the
two pieces $e_W$ and $e_E$, which may justify calling \emph{cut}
the corresponding inference rule.

We define now by induction the notion of \emph{construction of a global} K\emph{-graph}, which for short we call just \emph{construction}. A construction will be a finite binary tree in whose nodes we have an oriented graph together with some distinguished edges of this graph.

The oriented graph at the root of a construction $G$ will be called the \emph{root graph} of $G$, and we say that $G$ is a construction of its root graph. For $X\in\{W,E\}$ and $Y\in\{N,S\}$, we write $YX(G)$ for the distinguished edges of the root graph of $G$, which are at the root of $G$ together with the root graph.

Here are the two clauses of our definition of construction:
\begin{itemize}
\item[(1)]The single-node tree in whose single node, which is both the root and the unique leaf, we have the underlying oriented graph of a basic K-graph, together with the distinguished edges $YX(B)$, is a construction.
\item[(2)]For every $X\in\{W,E\}$, let $G_X$ be a construction of the oriented graph $D_X$, and let $e_X$ be a functional $\bar{X}$-edge of $D_X$. The tree of the construction $G=G_W[e_W\mn e_E]G_E$ is obtained by adding to the trees of the constructions $G_W$ and $G_E$ a new node, which will be the root of $G$, whose successors are the roots of the trees of $G_W$ and $G_E$. The oriented graph at the root of $G$, i.e.\ the root graph of $G$, is $D=D_W[e_W\mn e_E]D_E$ provided the following is satisfied for every $Y\in\{N,S\}$:\\[1.5ex]
(XYC)\hspace{3em}$e_X=Y\bar{X}(G_X)$\hspace{.5em} or \hspace{.5em}$e_{\bar{X}}=YX(G_{\bar{X}})$.\\[1.5ex]
(Note that this condition for $X$ being $W$ is the same as this condition for $X$ being $E$.) The distinguished edges of $D$ at the root of $G$ are obtained as follows:\\[1.5ex]
(XYD)\hspace{3em}$YX(G)=\left\{\begin{array}{ll}
YX(G_X) & {\mbox{\rm if }}\hspace{.5em} e_{\bar{X}}=YX(G_{\bar{X}}),
\\
YX(G_{\bar{X}}) & {\mbox{\rm otherwise.}}
\end{array}\right.$\\[1.5ex]
At the other nodes of the tree of $G$, which are not the root of $G$, we have in $G$ the same oriented graphs and
the same distinguished edges that we had in $G_X$.
\end{itemize}
This concludes our definition of construction.

A \emph{global} K\emph{-graph} is an oriented graph that is the root graph of a construction.

Let $G_W$ and $G_E$ be respectively constructions of the global
K-graphs $D_W$ and $D_E$ in the example for $D_W[e_W\mn e_E]D_E$
given above, and let $G=G_W[e_W\mn e_E]G_E$. We may take, as the
picture suggests, that $SE(G_W)=e_W$ and $NW(G_E)=e_E$. So (XYC)
would be satisfied. To illustrate how (XYD) is applied, we have
the following picture:
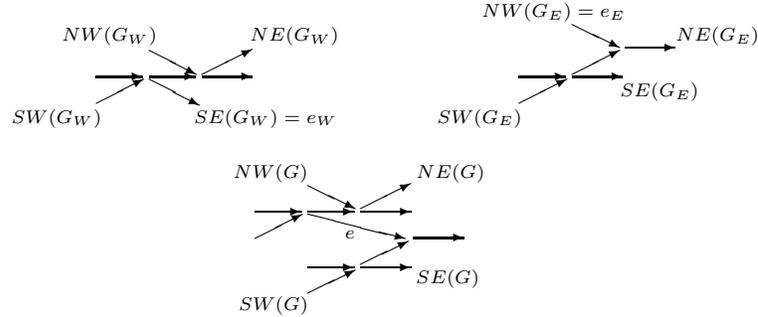
\begin{center}
\begin{picture}(200,50)
\put(0,10){\vector(2,1){18}} \put(0,20){\vector(1,0){17.5}}
\put(20,19){\vector(2,-1){19}} \put(20,20){\vector(1,0){17.5}}
\put(40,20){\vector(1,0){19}} \put(20,30){\vector(2,-1){18}}
\put(40,21){\vector(2,1){19}}

\put(64,8){\scriptsize\makebox(0,0)[t]{$SE(G_W)=e_W$}}
\put(-15,8){\scriptsize\makebox(0,0)[t]{$SW(G_W)$}}
\put(5,39){\scriptsize\makebox(0,0)[t]{$NW(G_W)$}}
\put(75,39){\scriptsize\makebox(0,0)[t]{$NE(G_W)$}}

\put(160,10){\vector(2,1){18}} \put(160,20){\vector(1,0){17.5}}
\put(180,21){\vector(2,1){18}} \put(180,20){\vector(1,0){19}}
\put(180,40){\vector(2,-1){18}} \put(200,31){\vector(1,0){19}}

\put(213,18){\scriptsize\makebox(0,0)[t]{$SE(G_E)$}}
\put(145,8){\scriptsize\makebox(0,0)[t]{$SW(G_E)$}}
\put(173,48){\scriptsize\makebox(0,0)[t]{$NW(G_E)=e_E$}}
\put(235,39){\scriptsize\makebox(0,0)[t]{$NE(G_E)$}}
\end{picture}
\end{center}
\begin{center}
\begin{picture}(80,60)
\put(0,30){\vector(2,1){18}} \put(0,40){\vector(1,0){17.5}}
\put(20,39){\vector(4,-1){37}} \put(20,40){\vector(1,0){17.5}}
\put(40,40){\vector(1,0){19}} \put(20,50){\vector(2,-1){18}}
\put(40,41){\vector(2,1){19}}

\put(20,9){\vector(2,1){18}} \put(20,19){\vector(1,0){17.5}}
\put(40,20){\vector(2,1){18}} \put(40,19){\vector(1,0){19}}
\put(60,30){\vector(1,0){19}}

\put(73,18){\scriptsize\makebox(0,0)[t]{$SE(G)$}}
\put(7,8){\scriptsize\makebox(0,0)[t]{$SW(G)$}}
\put(6,58){\scriptsize\makebox(0,0)[t]{$NW(G)$}}
\put(74,58){\scriptsize\makebox(0,0)[t]{$NE(G)$}}

\put(36,33){\scriptsize\makebox(0,0)[t]{$e$}}

\end{picture}
\end{center}

The terminology that we introduce now is based on \cite{H69}
(Chapter 16). For $n\geq 1$, consider a sequence $a_1,\dots,a_n$
of vertices of a digraph $D$ such that if $n\geq 2$, then for
every $i\in\{1,\ldots,n\mn 1\}$ we have that $(a_i,a_{i+1})\in D$
or $(a_{i+1},a_i)\in D$.

Such a sequence is a \emph{semipath} when all the vertices in it are mutually distinct, and it is a \emph{semicycle} when $a_1=a_n$, with $n\geq 4$, and all the vertices in $\{a_1\ldots,a_{n-1}\}$ are mutually distinct. If in the definition of semipath we replace the disjunction ``$(a_i,a_{i+1})\in D$ or $(a_{i+1},a_i)\in D$'' by the first disjunct ``$(a_i,a_{i+1})\in D$'', then we obtain the definition of \emph{path}.

When $a_1,\dots,a_n$ is a semipath, we say that $a_1$ is \emph{joined} to $a_n$ by the semipath $a_1,\dots,a_n$. Note that for every semipath $a_1,\dots,a_n$, the sequence in the inverse order $a_n,\dots,a_1$ is also a semipath. We call then $a_1,\dots,a_n$ and $a_n,\dots,a_1$ \emph{cognate semipaths}.

A digraph is \emph{weakly connected} when every two vertices in it
are joined by a semipath. A digraph is \emph{asemicyclic} when it
has no semicycles. A digraph is \emph{W-E-functional} when all its
$W$-edges and $E$-edges are functional.

It is straightforward to prove by induction on the number of inner
edges that every global K-graph is a weakly connected, asemicyclic
and $W$-$E$-functional oriented graph that has an inner vertex.

For $k_X$ being the number of $X$-edges of an arbitrary global K-graph $D$ at the root of a construction $G$, we can prove the following analogue of (XYB):
\begin{tabbing}
\hspace{1.5em}(XYG)\hspace{3em}if $k_X\geq 2$, then $NX(G)\neq SX(G)$.
\end{tabbing}
\noindent{\sc Proof of (XYG). } We proceed by induction on the number of inner edges of $D$. In the basis, when $G$ is a basic K-graph, we have (XYB). In the induction step we have three cases.

(1) If both $NX(G)$ and $SX(G)$ are from $D_X$, then
$e_{\bar{X}}=NX(G_{\bar{X}})=SX(G_{\bar{X}})$. By the induction
hypothesis we know that $D_{\bar{X}}$ has no other $X$-edge save
$e_{\bar{X}}$. So all the $X$-edges of $D$ are $X$-edges of $D_X$,
and then we apply the induction hypothesis to $G_X$.

(2) If both $NX(G)$ and $SX(G)$ are from $D_{\bar{X}}$, then we
apply the induction hypothesis to $G_{\bar{X}}$.

(3) If one of $NX(G)$ and $SX(G)$ is from $D_X$, while the other
is from $G_{\bar{X}}$, then (XYG) is trivial because $D_X$ and
$D_{\bar{X}}$ are disjoint digraphs.\qed

Our purpose next is to find conditions equivalent with (XYD) of clause (2) of the definition of construction above. These equivalent conditions will come handy for proofs later on. Note first that (XYD) amounts to the following two implications:
\begin{tabbing}
\hspace{1.5em}\=(XYD1)\hspace{3em}\=if $e_{\bar{X}}=YX(G_{\bar{X}})$, then $YX(G)=YX(G_X)$,\\[1ex]
\>(XYD2)\>if $YX(G)\neq YX(G_{\bar{X}})$, then $e_{\bar{X}}=YX(G_{\bar{X}})$.
\end{tabbing}
We infer easily the following from these two implications:
\begin{tabbing}
\hspace{1.5em}\=(XYD1)\hspace{3em}\=if $e_{\bar{X}}=YX(G_{\bar{X}})$, then $YX(G)=YX(G_X)$\kill

\>(D)\>$YX(G)=YX(G_X)$ or $YX(G)= YX(G_{\bar{X}})$.
\end{tabbing}
Then from (D) we infer easily for every $Z\in\{W,E\}$ that
\begin{tabbing}
\hspace{1.5em}\=(XYD1)\hspace{3em}\=if $e_{\bar{X}}=YX(G_{\bar{X}})$, then $YX(G)=YX(G_X)$\kill

\>(XY1)\>if $YX(G)\in D_Z$, then $YX(G)= YX(G_Z)$.
\end{tabbing}
From (XYD2) we also infer easily that
\begin{tabbing}
\hspace{1.5em}\=(XYD1)\hspace{3em}\=if $e_{\bar{X}}=YX(G_{\bar{X}})$, then $YX(G)=YX(G_X)$\kill

\>(XY2)\>if $YX(G)\in D_X$, then $e_{\bar{X}}=YX(G_{\bar{X}})$.
\end{tabbing}
So we have deduced (XY1) and (XY2) from (XYD).

We will now show that, conversely, we may deduce (XYD) from (XY1) and (XY2). Here is how we obtain (XYD2):
\begin{tabbing}
\hspace{1.5em}if $YX(G)\neq YX(G_{\bar{X}})$, \=then $YX(G)\notin D_{\bar{X}}$, by (XY1),\\*[.5ex]
\>then $YX(G)\in D_X$,\\*[.5ex]
\>then $e_{\bar{X}}=YX(G_{\bar{X}})$, by (XY2).
\end{tabbing}
We infer (D) from (XY1):
\begin{tabbing}
\hspace{1.5em}if $YX(G)\neq YX(G_X)$, \=then $YX(G)\notin D_X$, by (XY1),\\*[.5ex]
\>then $YX(G)\in D_{\bar{X}}$,\\*[.5ex]
\>then $YX(G)= YX(G_{\bar{X}})$, by (XY1),
\end{tabbing}
and we infer (XYD1) from (D):
\begin{tabbing}
\hspace{1.5em}if $e_{\bar{X}}=YX(G_{\bar{X}})$, \=then $YX(G)\neq
YX(G_{\bar{X}})$,\\*[.5ex]
\>then $YX(G)=YX(G_X)$, by (D).
\end{tabbing}
So (XY1) and (XY2) have the same force as (XYD).

The remainder of this section is devoted to proving for our notion of global K-graph a completeness result, which will help us for the equivalence proofs in later sections. We start first with three propositions that involve the equations that are assumed for planar polycategories (see the Introduction). The equations involved in Propositions 2.1 and 2.2 are like the equations of multicategories (see \cite{L89}, Section 3; analogous equations are also assumed for operads), while the equation involved in Proposition 2.3 is dual to that involved in Proposition 2.2.

Let $P$, $Q$ and $R$ be constructions, and let $e_W$ and $f_W$ be $E$-edges of the root graphs of $P$ and $Q$ respectively, while $e_E$ and $f_E$ are $W$-edges of the root graphs of $Q$ and $R$ respectively. Let $G_1$ be $(P[e_W\mn e_E]Q)[f_W\mn f_E]R$ and let $G_2$ be $P[e_W\mn e_E](Q[f_W\mn f_E]R)$. We can prove the following.

\prop{Proposition 2.1}{We have that $G_1$ is a construction iff $G_2$ is a construction. The root graphs of these constructions are the same and the distinguished edges of these root graphs at these roots are the same.}

\dkz In this proof we write $[-]$ for both ${[e_W\mn e_E]}$ and ${[f_W\mn f_E]}$, since it is clear from the context which we have in mind. We show first that if $(P[-]Q)[-]R$ is a construction, then $P[-](Q[-]R)$ is a construction.

We have, by (XYC), that
\begin{itemize}
\item[]$P[-]Q$ is a construction iff for every $Y\in\{N,S\}$ we have $e_W=YE(P)$ or $e_E=YW(Q)$,
\item[]$(P[-]Q)[-]R$ is a construction iff for every $Y\in\{N,S\}$ we have $f_W=YE(P[-]Q)$ or $f_E=YW(R)$.
\end{itemize}
We have that $f_W=YE(P[-]Q)$ and the fact that $f_W$ is in the root graph of $Q$ imply $f_W=YE(Q)$, by (XY1). Since $(P[-]Q)[-]R$ is a construction, we can conclude that $Q[-]R$ is a construction.

To show that $P[-](Q[-]R)$ is a construction it remains to verify that we have $e_W=YE(P)$ or $e_E=YW(Q[-]R)$. We have the implication
\begin{tabbing}
\hspace{1.5em}if $e_W\neq YE(P)$, then $e_E=YW(Q)$,
\end{tabbing}
since $P[-]Q$ is a construction. We also have
\begin{tabbing}
\hspace{1.5em}if $e_W\neq YE(P)$, \=then $YE(P[-]Q)=YE(P)$, by (XYD2),\\*[.5ex]
\>then $f_W\neq YE(P[-]Q)$, since $f_W$ is not an edge of the\\
\`root graph of $P$,\\[.5ex]
\>then $f_E=YW(R)$, since $(P[-]Q)[-]R$ is a construction,\\[.5ex]
\>then $YW(Q[-]R)=YW(Q)$, by (XYD1),\\[.5ex]
\>then $e_E=YW(Q[-]R)$,
\end{tabbing}
by the implication established above. So $P[-](Q[-]R)$ is a construction.

We proceed analogously to show that if $P[-](Q[-]R)$ is a construction, then $(P[-]Q)[-]R$ is a construction. It is clear that the root graphs of these two constructions are the same. It remains to establish that the distinguished edges of these root graphs at these roots are the same.

We have
\begin{tabbing}
\hspace{1.5em}\=$YE((P[-]Q)[-]R)=\left\{\begin{array}{ll}
YE(R) & {\mbox{\rm if }}\hspace{.5em} f_W=YE(P[-]Q),
\\
YE(P[-]Q) & {\mbox{\rm otherwise,}}
\end{array}\right.$\\[2ex]
\>$YE(P[-]Q)=\left\{\begin{array}{ll}
YE(Q) & {\mbox{\rm if }}\hspace{.5em} e_W=YE(P),
\\
YE(P) & {\mbox{\rm otherwise.}}
\end{array}\right.$
\end{tabbing}
Since $f_W$ is an edge of the root graph of $Q$, we have that
\begin{tabbing}
\hspace{1.5em}\=if $f_W=YE(P[-]Q)$, then $f_W=YE(Q)$, by (XY1),\\[.5ex]
\>if $f_W=YE(P[-]Q)$, then $e_W=YE(P)$, by (XY2),
\end{tabbing}
and we have that
\begin{tabbing}
\hspace{1.5em}\=if $f_W=YE(Q)$ and $e_W=YE(P)$, then $f_W=YE(P[-]Q)$,
\end{tabbing}
because if $e_W=YE(P)$, then $YE(Q)=YE(P[-]Q)$, as stated above. So we have
\begin{tabbing}
\hspace{1.5em}\=$YE((P[-]Q)[-]R)=\left\{\begin{array}{lll}
YE(R) & {\mbox{\rm if }}\hspace{.5em} f_W=YE(Q) & {\mbox{\rm and }}\hspace{.5em} e_W=YE(P),
\\
YE(Q) & {\mbox{\rm if }}\hspace{.5em} f_W\neq YE(Q) & {\mbox{\rm and }}\hspace{.5em} e_W=YE(P),
\\
YE(P) & {\mbox{\rm if }}\hspace{.5em} f_W\neq YE(Q) & {\mbox{\rm and }}\hspace{.5em} e_W\neq YE(P).
\end{array}\right.$
\end{tabbing}
On the other hand,
\begin{tabbing}
\hspace{1.5em}\=$YE(P[-](Q[-]R))=\left\{\begin{array}{ll}
YE(Q[-]R) & {\mbox{\rm if }}\hspace{.5em} e_W=YE(P),
\\
YE(P) & {\mbox{\rm otherwise,}}
\end{array}\right.$\\[2ex]
\>$YE(Q[-]R)=\left\{\begin{array}{ll}
YE(R) & {\mbox{\rm if }}\hspace{.5em} f_W=YE(Q),
\\
YE(Q) & {\mbox{\rm otherwise,}}
\end{array}\right.$
\end{tabbing}
which implies that
\[
YX((P[-]Q)[-]R)=YX(P[-](Q[-]R))
\]
when $X$ is $E$. We proceed in a dual manner when $X$ is $W$.\qed

Let $P$, $Q$ and $R$ be constructions, and let $e_W$ and $f_W$ be
different $E$-edges of the root graph of $P$, while $e_E$ and
$f_E$ are $W$-edges of the root graphs of $Q$ and $R$
respectively. Let $G_1$ be $(P[e_W\mn e_E]Q)[f_W\mn f_E]R$ and let
$G_2$ be $(P[f_W\mn f_E]R)[e_W\mn e_E]Q$.

{\sc Proposition 2.2} is formulated exactly as Proposition 2.1 for these new constructions $G_1$ and $G_2$.

\vspace{2ex}

\noindent{\sc Proof of Proposition 2.2. } As in the preceding
proof, we use the abbreviation $[-]$. We show that if
$(P[-]Q)[-]R$ is a construction, then $(P[-]R)[-]Q$ is a
construction. We have that $P[-]Q$ and $(P[-]Q)[-]R$ are
constructions under the same conditions concerning $e_X$ and $f_X$
displayed at the beginning of the proof of Proposition 2.1.

We have that $f_W=YE(P[-]Q)$ and the fact that $f_W$ is in the root graph of $P$ imply $f_W=YE(P)$, by (XY1). Since $(P[-]Q)[-]R$ is a construction, we can conclude that $P[-]R$ is a construction.

To show that $(P[-]R)[-]Q$ is a construction it remains to verify that we have $e_W=YE(P[-]R)$ or $e_E=YW(Q)$. We have
\begin{tabbing}
\hspace{1.5em}if $e_E\neq YW(Q)$, \=then $e_W=YE(P)$, since $P[-]Q$ is a construction,\\*[.5ex]
\>then $f_W\neq YE(P)$, since $e_W\neq f_W$,\\[.5ex]
\>then $YE(P[-]R)=YE(P)$, by (XYD2),\\[.5ex]
\>then $e_W=YE(P[-]R)$.
\end{tabbing}
So $(P[-]R)[-]Q$ is a construction.

We proceed in exactly the same manner to show that if $(P[-]R)[-]Q$ is a construction, then $(P[-]Q)[-]R$ is a construction. It is clear that the root graphs of these two constructions are the same. It remains to establish that the distinguished edges of these root graphs at these roots are the same.

We can conclude that
\begin{tabbing}
\hspace{.5em}\=$YE((P[-]Q)[-]R)=\left\{\begin{array}{lll}
YE(R) & {\mbox{\rm if }}\hspace{.5em} f_W=YE(P) & {\mbox{\rm and }}\hspace{.5em} e_W\neq YE(P),
\\
YE(Q) & {\mbox{\rm if }}\hspace{.5em} f_W\neq YE(P[-]Q) & {\mbox{\rm and }}\hspace{.5em} e_W=YE(P),
\\
YE(P) & {\mbox{\rm if }}\hspace{.5em} f_W\neq YE(P[-]Q) & {\mbox{\rm and }}\hspace{.5em} e_W\neq YE(P).
\end{array}\right.$
\end{tabbing}
Since $f_W$ is an edge of the root graph of $P$, we have that
\begin{tabbing}
\hspace{1.5em}if $f_W=YE(P[-]Q)$, then $f_W=YE(P)$ and $e_W\neq
YE(P)$,
\end{tabbing}
by using (XY1) and $e_W\neq f_W$. We also have that
\begin{tabbing}
\hspace{1.5em}if $f_W=YE(P)$ and $e_W\neq YE(P)$, then
$f_W=YE(P[-]Q)$,
\end{tabbing}
by the clause for $YE(P[-]Q)$.

We can conclude analogously that
\begin{tabbing}
\hspace{.5em}\=$YE((P[-]R)[-]Q)=\left\{\begin{array}{lll}
YE(Q) & {\mbox{\rm if }}\hspace{.5em} e_W=YE(P) & {\mbox{\rm and }}\hspace{.5em} f_W\neq YE(P),
\\
YE(R) & {\mbox{\rm if }}\hspace{.5em} e_W\neq YE(P[-]R) & {\mbox{\rm and }}\hspace{.5em} f_W=YE(P),
\\
YE(P) & {\mbox{\rm if }}\hspace{.5em} e_W\neq YE(P[-]R) & {\mbox{\rm and }}\hspace{.5em} f_W\neq YE(P).
\end{array}\right.$
\end{tabbing}

We show first that
\begin{tabbing}
\hspace{1.5em}\=(P)\hspace{3em}\=($f_W\neq YE(P[-]Q)$ and $e_W\neq YE(P)$) iff\\
\`($f_W\neq YE(P)$ and $e_W\neq YE(P[-]R)$).
\end{tabbing}
By (XY1), we have that
\begin{tabbing}
\hspace{1.5em}\=if $f_W=YE(P[-]Q)$, \=then $f_W=YE(P)$,\\[.5ex]
\>if $e_W=YE(P[-]R)$, \>then $e_W=YE(P)$,
\end{tabbing}
and the converse implications hold by the clauses for $YE(P[-]Q)$ and $YE(P[-]R)$ because $e_W\neq f_W$. This is enough to establish (P).

We establish that
\begin{tabbing}
\hspace{1.5em}\=(P)\hspace{3em}\=($f_W\neq YE(P[-]Q)$ and $e_W\neq YE(P)$) iff\kill

\>(Q)\>($f_W\neq YE(P[-]Q)$ and $e_W=YE(P)$) iff\\*
\`($f_W\neq YE(P)$ and $e_W=YE(P)$),\\[1ex]

\>(R)\>($f_W=YE(P)$ and $e_W\neq YE(P)$) iff\\*
\`($f_W=YE(P)$ and $e_W\neq YE(P[-]R)$),
\end{tabbing}
by using (XY1) and $e_W\neq f_W$. So we have that
\[
YE((P[-]Q)[-]R)=YE((P[-]R)[-]Q).
\]

We have that
\begin{tabbing}
\hspace{1.5em}\=$YW((P[-]Q)[-]R)=\left\{\begin{array}{lll}
YW(P) & {\mbox{\rm if }}\hspace{.5em} f_E=YW(R) & {\mbox{\rm and }}\hspace{.5em} e_E=YW(Q),
\\
YW(Q) & {\mbox{\rm if }}\hspace{.5em} f_E=YW(R) & {\mbox{\rm and }}\hspace{.5em} e_E\neq YW(Q),
\\
YW(R) & {\mbox{\rm if }}\hspace{.5em} f_E\neq YW(R),
\end{array}\right.$
\end{tabbing}
\begin{tabbing}
\hspace{1.5em}\=$YW((P[-]R)[-]Q)=\left\{\begin{array}{lll}
YW(P) & {\mbox{\rm if }}\hspace{.5em} f_E=YW(R) & {\mbox{\rm and }}\hspace{.5em} e_E=YW(Q),
\\
YW(R) & {\mbox{\rm if }}\hspace{.5em} e_E=YW(Q) & {\mbox{\rm and }}\hspace{.5em} f_E\neq YW(R),
\\
YW(Q) & {\mbox{\rm if }}\hspace{.5em} e_E\neq YW(Q).
\end{array}\right.$
\end{tabbing}
We have that
\begin{tabbing}
\hspace{1.5em}if $e_E\neq YW(Q)$, \=then $e_W=YE(P)$, \=since $P[-]Q$ is a construction,\\*[.5ex]
\>then $f_W\neq YE(P)$, \>since $e_W\neq f_W$,\\[.5ex]
\>then $f_E=YW(R)$, \>since $P[-]R$ is a construction,
\end{tabbing}
and by contraposition we have that if $f_E\neq YW(R)$, then $e_E=YW(Q)$. This, together with what we have established previously, shows that
\[
YX((P[-]Q)[-]R)=YX((P[-]R)[-]Q)
\]
for every $X\in\{W,E\}$ and every $Y\in\{N,S\}$.\qed

Let $P$, $Q$ and $R$ be constructions, and let $e_E$ and $f_E$ be different $W$-edges of the root graph of $P$, while $e_W$ and $f_W$ are $E$-edges of the root graphs of $Q$ and $R$ respectively. Let $G_1$ be $R[f_W\mn f_E](Q[e_W\mn e_E]P)$ and let $G_2$ be $Q[e_W\mn e_E](R[f_W\mn f_E]P)$.

{\sc Proposition 2.3} is formulated exactly as Proposition 2.1 for these new constructions $G_1$ and $G_2$. It is proved in a manner dual to what we had for the proof of Proposition 2.2.

Consider the relations between constructions that exist between the constructions $G_1$ and $G_2$ of Propositions 2.1, 2.2 and 2.3. We call these relations $\rho_1$, $\rho_2$ and $\rho_3$ respectively.

Let $\rho$\emph{-equivalence} be the equivalence relation between constructions that is the reflexive, symmetric and transitive closure of $\rho_1\cup\rho_2\cup\rho_3$, and which is closed moreover under $\rho$\emph{-congruence}:
\begin{itemize}
\item[]if $G_1$ is $\rho$-equivalent with $G_2$ and $H_1$ is $\rho$-equivalent with $H_2$, then\\ $G_1[e_W\mn e_E]H_1$ is $\rho$-equivalent with $G_2[e_W\mn e_E]H_2$,
\end{itemize}
provided the last two constructions are defined. We can prove the following.

\prop{Proposition 2.4}{If $e$ is an inner edge of the root graph of a construction $G$, then there are two constructions $H_W$ and $H_E$ such that $G$ is $\rho$-equivalent to $H_W[e_W\mn e_E]H_E$.}

\dkz We proceed by induction on the number $n$ of inner edges in
the root graph of $G$. If $n=1$, then $G$ is of the form
$G_W[e_W\mn e_E]G_E$, and we take $H_X$ to be $G_X$. If $n\geq 2$,
and $G$ is again of that form, then again we choose $H_X$ to be
$G_X$.

Suppose $n\geq 2$, and $G$ is of the form $G_W[f_W\mn f_E]G_E$ for
$f$ different from $e$. If $e$ is in the root graph of $G_W$, then
by the induction hypothesis $G$ is $\rho$-equivalent to a
construction
\[
(G_{WW}[e_W\mn e_E]G_{WE})[f_W\mn f_E]G_E,
\]
which is $\rho$-equivalent to either
\[
G_{WW}[e_W\mn e_E](G_{WE}[f_W\mn f_E]G_E),
\]
because of $\rho_1$, or
\[
(G_{WW}[f_W\mn f_E]G_E)[e_W\mn e_E]G_{WE},
\]
because of $\rho_2$. We proceed analogously if $e$ is in the root graph of $G_E$, by appealing to $\rho_1$ and $\rho_3$. \qed

We say that a basic K-graph $B$ \emph{determines} a leaf of a construction $G$ when $B$ occurs in an application of clause (1) for the definition of $G$. For a given construction $G$, let $[G]$ be the set of all the constructions that have leaves determined by the same basic K-graphs as $G$, and that have the same root graph as $G$. We can prove the following for every pair of constructions $G$ and $H$.

\prop{Proposition 2.5}{We have that $G$ and $H$ are $\rho$-equivalent iff $[G]=[H]$.}

\dkz For the implication from left to right we have essentially just an easy application of Propositions 2.1-2.3. For the other direction, suppose $[G]=[H]$. We proceed by induction on the number $n$ of inner edges in the root graph $D$ of $G$ and $H$, which they share. If $n=0$, then $G$ and $H$ are the same construction, given by the same basic K-graph.

Let $n\geq 1$, and let $e$ be an inner edge of $D$. Then by Proposition 2.4 we have that $G$ and $H$ are $\rho$-equivalent to respectively $G_W[e_W\mn e_E]G_E$ and $H_W[e_W\mn e_E]H_E$. We apply the induction hypothesis to $G_X$ and $H_X$, and then we appeal to $\rho$-congruence. \qed

With that we have proved the completeness result we set ourselves as a goal in this section. As a consequence of Propositions 2.1-2.3 we also have the following for every pair of constructions $G$ and $H$, for every $X\in\{W,E\}$ and every $Y\in\{N,S\}$.

\prop{Proposition 2.6}{If $G$ and $H$ are $\rho$-equivalent, then $YX(G)=YX(H)$.}

We conclude this section with some terminological matters, which we need for the exposition later on. In a construction $G$ let the \emph{root vertices} of $G$ be the vertices of the root graph of $G$. The other vertices that may occur in the oriented graph at a node of $G$ that is not the root, which are not root vertices, will be called \emph{secondary vertices}.

Two constructions are said to be $\sigma$\emph{-equivalent} when they are in all respects the same, save that they may differ in the choice of secondary vertices. One could say that they are the same construction up to renaming of secondary vertices.

For $G$ a construction, consider $[G]$, and let $\|G\|$ be the set of all the constructions $\sigma$-equivalent to a construction in $[G]$. We call $\|G\|$ a \emph{global compass graph}.

\section{Local K-graphs}
In this section we deal with our second definition of K-graph (see
the Introduction), which yields the notion of local K-graph.

For $D$ a digraph and $a$ an inner vertex of $D$ consider for
$Y\in\{N,S\}$ the two functions $YW$ such that $YW(a)$ is an edge
of $D$ of the form $(b,a)$, and consider the two functions $YE$
such that $YE(a)$ is an edge of $D$ of the form $(a,b)$. For every
inner vertex $a$ of $D$ let $k^a_W\geq 1$ be the number of edges
of $D$ of the form $(b,a)$, while $k^a_E\geq 1$ is the number of
edges of $D$ of the form $(a,b)$.

Let $\L$ be a set of such four functions. Then we say that
$\langle D,\L\rangle$ \emph{separates N from S} when the following
condition (analogous to (XYB) of Section~2) holds for
every inner vertex $a$ of $D$ and every $X\in\{W,E\}$:
\begin{tabbing}
\hspace{1.5em}if $k^a_X\geq 2$, then $NX(a)\neq SX(a)$.
\end{tabbing}

We say that a path $a_1,\ldots,a_n$, with $n\geq 1$, of $D$ is
$Y$-\emph{decent} in $\langle D,\L\rangle$ when either $n=1$ or if
$n\geq 2$, then $YE(a_1)=(a_1,a_2)$ or $YW(a_n)=(a_{n-1},a_n)$
(see Section~2 for the notion of path). A path of $D$ is
\emph{decent} in $\langle D,\L\rangle$ when it is both $N$-decent
and $S$-decent in $\langle D,\L\rangle$.

For example, if $\langle D,\L\rangle$ is such that $NE(a)=(a,d)$
and $NW(c)=(e,c)$, as in the following picture of $D$, then the
path $a,b,c$ is not $N$-decent in $\langle D,\L\rangle$:
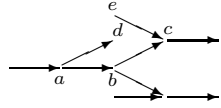
\begin{center}
\begin{picture}(60,50)
\put(0,20){\vector(1,0){19}} \put(20,20){\vector(1,0){19}}
\put(20,21){\vector(2,1){19}} \put(40,9){\vector(1,0){17.5}}
\put(60,9){\vector(1,0){19}} \put(40,19){\vector(2,-1){18}}
\put(40,21){\vector(2,1){18}} \put(40,40){\vector(2,-1){18}}
\put(60,30.5){\vector(1,0){19}}

\put(19,14){\scriptsize\makebox(0,0)[b]{$a$}}
\put(39,13.5){\scriptsize\makebox(0,0)[b]{$b$}}
\put(60,33){\scriptsize\makebox(0,0)[b]{$c$}}
\put(41,37){\scriptsize\makebox(0,0)[t]{$d$}}
\put(39,45){\scriptsize\makebox(0,0)[t]{$e$}}

\end{picture}
\end{center}

For $D$ an oriented graph, we say that $\langle D,\L\rangle$ is a
\emph{local compass graph} when
\begin{tabbing}
\hspace{1.5em}\=(XYb)\hspace{1em}\=if \kill

\>(1) \>$D$ is weakly connected,
\\[.5ex]
\>(2) \>$D$ is asemicyclic,
\\[.5ex]
\>(3) \>$D$ is $W$-$E$-functional and has an inner vertex,
\\[.5ex]
\>(4) \>$\langle D,\L\rangle$ separates $N$ from $S$,
\\[.5ex]
\>(5) \>every path of $D$ is decent in $\langle D,\L\rangle$.
\end{tabbing}
For a local compass graph $\langle D,\L\rangle$ we say that the
oriented graph $D$ is a \emph{local} K-\emph{graph}.

We say that a path $a_1,\ldots,a_n$, with $n\geq 2$, of a digraph
$D$ \emph{covers} an edge $e$ of $D$ when $e=(a_i,a_{i+1})$ for
some $i\in\{1,\ldots,n\mn 1\}$. We need this notion for the
inductive definition of local compass graph, which we will now give.

It is clear that a basic K-graph $B$ with the unique inner vertex
$b$ gives rise to a local compass graph $\langle D,\L\rangle$
where $D$ is the oriented graph underlying $B$ and $YX(b)=YX(B)$.
Starting from these local compass graphs we could define local
compass graphs inductively.

If $\langle D_W,\L_W\rangle$ and $\langle D_E,\L_E\rangle$ are
local compass graphs, then for $D=D_W[e_W\mn e_E]D_E$ we have that
$\langle D,\L\rangle$ is a local compass graph provided the
oriented graph $D=D_W[e_W\mn e_E]D_E$ is defined, and the
functions in $\L$ are defined by taking that for an inner vertex
$a$ of $D_Z$, where $Z\in\{W,E\}$, we have that $YX(a)$ has the
same value as in $\langle D_Z,\L_Z\rangle$ if this value is
different from $e_W$ and $e_E$; otherwise it is $e$. We assume
moreover a condition that will yield (5) above:
\begin{tabbing}
\hspace{1.5em}\=every path of $D$ that covers $e$ is decent in
$\langle D,\L\rangle$.
\end{tabbing}
The conditions (1)-(4) are then easily derived.

The equivalence of the two notions of local compass graph, the one
given by the first definition, in terms of (1)-(5), and the one given by the second,
inductive, definition is established in a straightforward manner.

\section{Global and local K-graphs}
In this section we establish the equivalence between the notions
of global and local K-graph.

For a given construction $G$ we define the local compass graph
$\lambda(G)=\langle D,\L\rangle$ in the following manner. The
oriented graph $D$ is the root graph of $G$, and the functions in
$\L$ are defined inductively. In the basis, for a basic K-graph
$B$ with the inner vertex $b$ we have that $YX(b)$ is defined as
$YX(B)$. In the induction step, if $G$ is $G_W[e_W\mn e_E]G_E$ and, for $Z\in\{W,E\}$, we have $\lambda(G_Z)=\langle
D_Z,\L_Z\rangle$, then $YX(a)$ in $\L$ has the same value as
$YX(a)$ in $\L_Z$, provided $a$ is in $D_Z$, except when this
value was $e_W$ or $e_E$, in which case the value is $e$ in~$\L$.

This inductive definition of $\L$ makes the functions $YX$ in it
dependent only on the arrangement of distinguished edges of the
basic K-graphs in the leaves of $G$. Hence $\lambda(G)$ depends
only on this arrangement and on the root graph of~$G$.

So we could define a function $\Lambda$ from global compass graphs
$\|G\|$ (see the end of Section~2 for $\|G\|$ and $[G]$) to local compass graphs such
that $\Lambda\|G\|=\lambda(G)$. It is easy to verify that if
$\|G\|=\|H\|$, then $\lambda(G)=\lambda(H)$, which implies that if
$[G]=[H]$, then $\lambda(G)=\lambda(H)$. It remains to verify that
$\lambda(G)$ is indeed a local compass graph.

For this verification, we have established (1)-(3) of the definition of local
compass graph in Section~2, and condition (4) of this definition
is immediate from (XYB). It remains to verify condition (5). For
that we need some preliminary matters. The following definitions
apply to $\lambda(G)=\langle D,\L\rangle$ as defined above, but
the same definitions may be given for every local compass graph
$\langle D,\L\rangle$.

We say in $\langle D,\L\rangle$ that $(a,b)$ is a $YW$-\emph{edge}
of $D$ when $YW(b)=(a,b)$ or $(a,b)$ is an $E$-edge of $D$, and we
say that $(a,b)$ is a $YE$-\emph{edge} of $D$ when $YE(a)=(a,b)$
or $(a,b)$ is a $W$-edge of $D$. (The $E$-edges and $W$-edges of
$D$ are both functional.) A path of $D$ is a $YX$-\emph{path} when
every edge it covers is a $YX$-edge.

We can prove the following for every construction $G$, for every
$Y\in\{N,S\}$ and every $X\in\{W,E\}$.

\begin{tabbing}
\textsc{Proposition
4.1.}\hspace{.25em}\=($W$)\hspace{.25em}\=\emph{If}\hspace{.25em}\=
$YW(G)\:$\=$=(a,b)$ \emph{in} $G$, \=\emph{then}
$YW(b)\:$\=$=(a,b)$ \emph{in} $\lambda(G)$.
\\*[.5ex]
\>($E$)\>\emph{If}\>$YE(G)\;$\>$=(a,b)$ \emph{in} $G$, \>\emph{then}
$YE(a)$\>$=(a,b)$ \emph{in} $\lambda(G)$.
\end{tabbing}

\dkz For ($W$), we proceed by induction on the number $n$ of inner
edges of the root graph $D$ of $G$. In the basis, when $n=0$, we
deal with a basic K-graph, and the implication holds trivially. If
$n\geq 1$, and $e$ is an inner edge of $D$, then by Proposition
2.4 we have that $G$ is $\rho$-equivalent to $H=H_W[e_W\mn
e_E]H_E$, and by Proposition 2.6, we have that $YX(G)=YX(H)$. For
$D_Z$ being the root graph of $H_Z$, we have

\begin{tabbing}
\hspace{1.5em}if $(a,b)=YW(H)\in D_Z$, \=then
$YW(H)=YW(H_Z)$, by (XY1),
\\[.5ex]
\>then $YW(b)=(a,b)$ in $\lambda(H_Z)$,
\end{tabbing}
by the induction hypothesis applied to $H_Z$. It is then clear that $YW(b)=(a,b)$ in $\lambda(G)$, since the
basic K-graphs of $H_Z$ are taken over by $G$. We prove ($E$)
analogously \qed

\oprop{Proposition 4.2}{We have that $h=YX(G)$ iff $h$ is an
$X$-edge of the root graph $D$ of $G$ such that every path of $D$
that covers $h$ is a $YX$-path in $\lambda(G)$.}

\dkz Suppose $X$ is $W$. From left to right we proceed by
induction on the number $n$ of inner edges of $D$. In the basis,
when $n=0$, we deal with a basic K-graph, and the proposition
holds trivially. If $n\geq 1$, consider a path $a_1,\ldots,a_m$,
with $m\geq 2$, that covers $YW(G)$. If for $(a_i,a_{i+1})=e$,
where $i\in\{1,\ldots,m\mn 1\}$, we have that
$YW(a_{i+1})=(c,a_{i+1})\neq e$, by Proposition 2.4 we have that
$G$ is $\rho$-equivalent to ${H=H_W[e_W\mn e_E]H_E}$. (Note that
$e$ must be an inner edge of~$D$, by (W) of Proposition 4.1.)

Since $e_E\neq YW(a_{i+1})$, we have $e_E\neq YW(H_E)$, by ($W$)
of Proposition 4.1. So $YW(H)=YW(H_E)$, by (XYD), which, together
with Proposition 2.6, contradicts the assumption that $h=YW(G)$.

From right to left we make again an induction on the number $n$ of
inner edges of $D$. The basis, when $n=0$, is again trivial. For
the induction step, when $n\geq 1$, suppose $h\neq YW(G)$. We want
to show that if $h$ is a $W$-edge of $D$, then there is a path
$a_1,\ldots,a_m$, with $m\geq 2$, such that $(a_1,a_2)=h$ and this
path is not a $YW$-path. Since $D$ is weakly connected, there is a
semipath $a_1, \ldots, a_m, b_2, \ldots, b_k$ of $D$, with $m\geq
2$, $a_m=b_1$ and $k\geq 2$, such that $(a_1,a_2)=h$,
$(b_k,b_{k-1})=YW(G)$, $a_1,\ldots,a_m$ is a path of $D$ and
$(b_2,b_1)=e\in D$. We use here the assumption that $h$ is a
$W$-edge of $D$; otherwise, $YW(G)$ could, for example, be of the
form $(c,a_1)$.

If $e=YW(G)$, then, by ($W$) of Proposition 4.1, we have
$e=YW(a_m)$ in $\lambda(G)$, and $(a_{m-1},a_m)$ is not a
$YW$-edge. If $e\neq YW(G)$, then, with the help of the assumption
that $D$ is $W$-$E$-functional, we conclude that $e$ is an inner
edge of $D$, and then, by Proposition 2.4, we have that $G$ is
$\rho$-equivalent to $H=H_W[e_W\mn e_E]H_E$. We must have that
$YW(G)$, which, by Proposition 2.6, is equal to $YW(H)$, is in the
root graph of $H_W$ (because we have a semipath $b_2,\ldots,b_k$
in this root graph). By (XY2), we conclude that $e_E=YW(H_E)$,
and, by ($W$) of Proposition 4.1, we have
$e_E=YW(a_m)\neq(a_{m-1},a_m)$. So $a_1,\ldots,a_m$ is not a
$YW$-path in $\lambda(H_E)$, which implies that it is not a
$YW$-path in $\lambda(G)$. We proceed analogously when $X$ is $E$.
\qed

We can now prove the following for every construction $G$.

\prop{Proposition 4.3}{We have that $\lambda(G)$ is a local
compass graph.}

\dkz As we noted at the beginning of the section, it remains to
verify condition (5) of the definition of local compass graph.

Suppose we have a path $a_1,\ldots,a_n$ of the root graph $D$ of
$G$ that is not decent. So $n\geq 2$, and for some $Y\in\{N,S\}$
we have $YE(a_1)\neq(a_1,a_2)$ and $YW(a_n)\neq(a_{n-1},a_n)$. Any
edge covered by this path must be an inner edge of $D$, and, since
$n\geq 2$, there is such an edge; let us call it $e$. By
Proposition 2.4, we have that $G$ is $\rho$-equivalent to
$H_W[e_W\mn e_E]H_E$. By Proposition 4.2, we conclude that
$e_W\neq YE(H_W)$ and $e_E\neq YW(H_E)$, but this contradicts the
fact that $H_W[e_W\mn e_E]H_E$ is a construction. \qed

The following two propositions serve to prove that there is a
bijection between global and local compass graphs

\prop{Proposition 4.4}{For every local compass graph $\langle
D,\L\rangle$ there is a construction $G$ such that
$\lambda(G)=\langle D,\L\rangle$.}

\dkz We proceed by induction on the number $n$ of inner edges of
$D$. If $n=0$, then $\langle D,\L\rangle$ determines a basic
K-graph, and the proposition holds trivially. If $n\geq 1$, then
$D$ is of the form $D_W[e_W\mn e_E]D_E$ for the local compass
graphs $\langle D_W,\L_W\rangle$ and $\langle D_E,\L_E\rangle$.
This follows from the inductive definition of local compass
graphs, which gives an equivalent notion.

By the induction hypothesis, for $X\in\{W,E\}$ we have the
constructions $G_X$ such that $\lambda(G_X)=\langle
D_X,\L_X\rangle$. We show first that $G=G_W[e_W\mn e_E]G_E$ is a
construction of a global K-graph. For that we have to check (XYC).
Suppose for some $Y\in\{N,S\}$ we have $e_W\neq YE(G_W)$ and
$e_E\neq YW(G_E)$. By Proposition 4.2, there is a path
$a_1,\ldots,a_n$, where $n\geq 2$, in $D_W$ that is not a $YE$-path
with $(a_{n-1},a_n)=e_W$, and there is a path $b_1,\ldots,b_m$,
where $m\geq 2$, in $D_E$ that is not a $YW$-path with
$(b_1,b_2)=e_E$. We may assume that $(a_1,a_2)\neq YE(a_1)$ and
$(b_{m-1},b_m)\neq YW(b_m)$. The path
$a_1,\ldots,a_{n-1},b_2,\ldots,b_m$ of $D$ is not a decent path.
So (XYC) holds.

To finish the proof we have to check that $\lambda(G)=\langle
D,\L\rangle$. It is clear that the root graph of $G_W[e_W\mn
e_E]G_E$ is $D$, while the definition of $\L$ in terms of $\L_W$
and $\L_E$ involved in the definition of $\lambda(G)$ is in
accordance with the clause for $\L$ in the inductive definition of
local compass graph. \qed

\oprop{Proposition 4.5}{If $\lambda(G)=\lambda(H)$, then
$\|G\|=\|H\|$.}

\dkz Let $\lambda(G)=\lambda(H)=\langle D,\L\rangle$. The oriented
graph $D$ together with the functions in $\L$ determines the basic
K-graphs that enter into the inductive definitions of $G$ and $H$
up to renaming of secondary vertices (see the end of Section~2).
Since $D$ is the root graph of both $G$ and $H$, we may conclude
that $[G]$ and $[H]$ are the same up to renaming of these
secondary vertices, which means that $\|G\|=\|H\|$. \qed

If we define $\Lambda\|G\|$ as $\lambda(G)$, as we did at the
beginning of this section, then from Propositions 4.4 and 4.5 we
infer that $\Lambda$ is a bijection between global and local
compass graphs. From the definition of this bijection, we may
conclude that the notions of global and local K-graphs coincide.

\section{K-graphs}
In this section we deal with our third definition of K-graph (see
the Introduction). For the notion this definition gives we
establish that it is equivalent with the notion of local K-graph,
and hence, by the results of Section~4, with both notions given by the preceding two
definitions.

The following definitions are for oriented graphs, and build upon
notions defined in Section~2. A \emph{proper} semipath is a
semipath such that neither it nor its cognate is a path.
Intuitively, there must be a change of direction in a proper semipath.

An edge $(a,b)$ is \emph{transversal} when there is a proper
semipath $a,b,\ldots,c$ and a proper semipath $b,a,\ldots,d$.

A \emph{bifurcation} is a triple of different edges that have a
common vertex. The following four kinds of bifurcations are
possible:
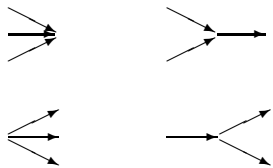
\begin{center}
\begin{picture}(100,35)
\put(0,10){\vector(2,1){18}} \put(0,30){\vector(2,-1){18}}
\put(0,20){\vector(1,0){17.5}}

\put(60,10){\vector(2,1){18}} \put(60,30){\vector(2,-1){18}}
\put(79,20){\vector(1,0){18}}
\end{picture}
\end{center}

\begin{center}
\begin{picture}(100,28)
\put(0,21){\vector(2,1){19}} \put(0,19){\vector(2,-1){19}}
\put(0,20){\vector(1,0){19}}

\put(80,21){\vector(2,1){19}} \put(80,19){\vector(2,-1){19}}
\put(60,20){\vector(1,0){19}}
\end{picture}
\end{center}
A bifurcation is called \emph{transversal} when all the three
edges in it are transversal.

A K-\emph{graph} is an oriented graph $D$ such that we have (1),
(2) and (3) from the definition of local K-graph of Section~3, and
we have moreover (instead of (4) and (5)) the following condition:
\begin{tabbing}
\centerline{No bifurcation is transversal.}
\end{tabbing}

In a semipath $a_1,\ldots,a_n$, with $n\geq 2$, of an oriented
graph $D$ we have for $i\in\{1,\ldots,n\mn 1\}$ that either
$(a_i,a_{i+1})$ or $(a_{i+1},a_i)$ is an edge of $D$, but not
both. We call this edge of $D$ the edge that \emph{connects} $a_i$
and $a_{i+1}$. We can now prove the following.

\prop{Proposition 5.1}{If in a semipath $a_1,\ldots, a_n$, with
$n\geq 2$, of an asemicyclic oriented graph $D$ the edge that
connects $a_1$ and $a_2$ and the edge that connects $a_{n-1}$ and
$a_n$ are transversal, then for every $i\in\{1,\ldots,n\mn 1\}$
the edge that connects $a_i$ and $a_{i+1}$ is transversal.}

\dkz If the edge that connects $a_1$ and $a_2$ is transversal,
then there is a proper semipath $a_2,a_1,\ldots,c$, and if the
edge that connects $a_{n-1}$ and $a_n$ is transversal, then there
is a proper semipath $a_{n-1},a_n,\ldots,d$. For every
$i\in\{1,\ldots,n\mn 1\}$ we have that
\[
\begin{array}{l}
a_{i+1},a_i,\ldots,a_2,a_1,\ldots,c
\\[.5ex]
a_i,a_{i+1},\ldots,a_{n-1},a_n,\ldots,d
\end{array}
\]
are proper semipaths. They are semipaths because $D$ is
asemicyclic, and hence all their members are mutually distinct,
and they are proper because they extend proper semipaths. We can
conclude that the edge that connects $a_i$ and $a_{i+1}$ is
transversal. \qed

For every K-graph $D$, if $D$ has transversal edges, by relying on
Proposition 5.1, we conclude that all the transversal edges of $D$
make a semipath $a_1,\ldots,a_n$, with $n\geq 2$, which we will
call the \emph{transversal} of $D$. The transversal is unique up
to cognation; the transversal is either a semipath or its cognate
(see Section~2). The vertices in the transversal, which must all
be inner, are called \emph{transversal vertices}.

All the edges of $D$ that share a single vertex with the transversal of $D$
are either \emph{in-going}, when for some $i\in\{1,\ldots,n\}$
they are of the form $(b,a_i)$, or they are \emph{out-going}, when
they are of the form $(a_i,b)$, where $a_i$ is a transversal
vertex. For an in-going edge $(b,a_i)$ we have in $D$ a tree
oriented from the leafs towards the root $a_i$:
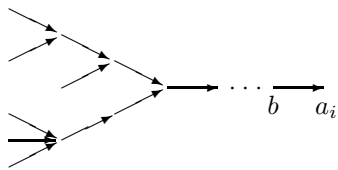
\begin{center}
\begin{picture}(120,77)
\put(0,10){\vector(2,1){18}} \put(0,30){\vector(2,-1){18}}
\put(0,20){\vector(1,0){17.5}}

\put(20,20){\vector(2,1){19}} \put(40,30){\vector(2,1){18}}

\put(0,70){\vector(2,-1){18}} \put(20,60){\vector(2,-1){18}}
\put(40,50){\vector(2,-1){18}} \put(0,50){\vector(2,1){18}}
\put(20,40){\vector(2,1){18}}

\put(60,40){\vector(1,0){19}} \put(100,40){\vector(1,0){19}}

\multiput(85,40)(5,0){3}{\makebox(0,0){\circle*{1}}}

\put(100,30){\makebox(0,0)[b]{$b$}}
\put(120,30){\makebox(0,0)[b]{$a_i$}}
\end{picture}
\end{center}
which we call an \emph{in-going tree}.

For an out-going edge $(a_i,b)$ we have in $D$ a tree oriented
from the root $a_i$ towards the leafs, which we call an
\emph{out-going tree}. These trees cannot share an edge with the
transversal of $D$; all the vertices in these trees except $a_i$
are not in the transversal of $D$. The orientation is imposed
because no transversal edge of $D$ is in these trees. If our
K-graph does not have transversal edges, then it has no
transversal, and is made only of trees analogous to in-going and
out-going trees that share a root.

The following proposition establishes that the notions of K-graph
and local K-graph are equivalent.

\prop{Proposition 5.2}{An oriented graph is a local K-graph iff it
is a K-graph.}

\dkz Let $D$ be an oriented graph that satisfies (1), (2) and (3)
of Section~3. To prove the proposition from left to right, suppose
that there is a transversal bifurcation in $D$. This bifurcation
can be of the four kinds mentioned above, which will produce in
$D$ subgraphs of the following four patterns (a subgraph of a
digraph is given by a subset of its edges on a subset of its
vertices):
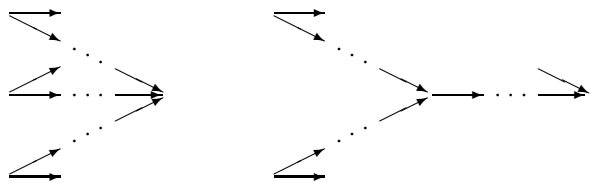
\begin{center}
\begin{picture}(220,80)
\put(0,10){\vector(2,1){19}} \put(0,9){\vector(1,0){19}}
\put(0,41){\vector(2,1){19}} \put(0,40){\vector(1,0){19}}
\put(0,70){\vector(2,-1){19}} \put(0,71){\vector(1,0){19}}

\multiput(25,22.5)(5,2.5){3}{\makebox(0,0){\circle*{1}}}
\multiput(25,57.5)(5,-2.5){3}{\makebox(0,0){\circle*{1}}}
\multiput(25,40)(5,0){3}{\makebox(0,0){\circle*{1}}}

\put(40,30){\vector(2,1){18}} \put(40,40){\vector(1,0){17.5}}
\put(40,50){\vector(2,-1){18}}

\put(100,10){\vector(2,1){19}} \put(100,9){\vector(1,0){19}}
\put(200,50){\vector(2,-1){19}} \put(160,40){\vector(1,0){19}}
\put(100,70){\vector(2,-1){19}} \put(100,71){\vector(1,0){19}}

\multiput(125,22.5)(5,2.5){3}{\makebox(0,0){\circle*{1}}}
\multiput(125,57.5)(5,-2.5){3}{\makebox(0,0){\circle*{1}}}
\multiput(185,40)(5,0){3}{\makebox(0,0){\circle*{1}}}

\put(140,30){\vector(2,1){18}} \put(200,40){\vector(1,0){17.5}}
\put(140,50){\vector(2,-1){18}}
\end{picture}
\end{center}
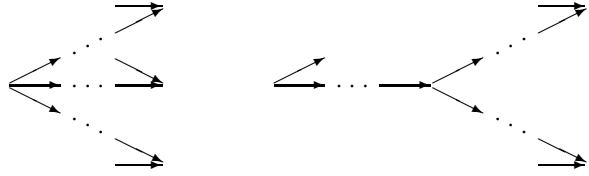
\begin{center}
\begin{picture}(220,80)
\put(40,20){\vector(2,-1){18}} \put(40,10){\vector(1,0){17.5}}
\put(40,50){\vector(2,-1){18}} \put(40,40){\vector(1,0){17.5}}
\put(40,60){\vector(2,1){18}} \put(40,70){\vector(1,0){17.5}}

\multiput(25,27.5)(5,-2.5){3}{\makebox(0,0){\circle*{1}}}
\multiput(25,52.5)(5,2.5){3}{\makebox(0,0){\circle*{1}}}
\multiput(25,40)(5,0){3}{\makebox(0,0){\circle*{1}}}

\put(0,41){\vector(2,1){19}} \put(0,40){\vector(1,0){19}}
\put(0,39){\vector(2,-1){19}}

\put(200,20){\vector(2,-1){18}} \put(200,10){\vector(1,0){17.5}}
\put(200,60){\vector(2,1){18}} \put(200,70){\vector(1,0){17.5}}

\multiput(185,27.5)(5,-2.5){3}{\makebox(0,0){\circle*{1}}}
\multiput(185,52.5)(5,2.5){3}{\makebox(0,0){\circle*{1}}}
\multiput(125,40)(5,0){3}{\makebox(0,0){\circle*{1}}}

\put(160,41){\vector(2,1){19}} \put(160,39){\vector(2,-1){19}}

\put(100,40){\vector(1,0){19}} \put(140,40){\vector(1,0){19}}
\put(100,41){\vector(2,1){19}}

\end{picture}
\end{center}
(These four subgraphs play here a role analogous to Kuratowski's
graphs $K_5$ and $K_{3,3}$, one of which must be found in
nonplanar graphs; see \cite{H69}, Chapter 11. Actually, if in the
graphs where these subgraphs occur there is a single $W$-vertex
and a single $E$-vertex, then in these graphs we have an oriented
version of $K_{3,3}$. On the other hand, $K_5$ is related to
asemicyclicity. We intend to deal with these matters one another
occasion.)

In all the four cases we go through all possible functions that
could make $\L$ to show that there must be a path of $D$ that is
not decent in $\langle D,\L\rangle$. This establishes the
proposition from left to right.

To prove the proposition from right to left, assume we are given a K-graph $D$. We define the functions in $\L$ by giving their value first for
non-transversal inner vertices $b$. We can do it in many ways,
provided we take care to guarantee that $\langle D,\L\rangle$
separates $N$ from $S$ (see Section~3). For $X\in\{W,E\}$ we choose
$NX(b)$ and $SX(b)$ as the same edge when there is no $X$-ward
branching in $b$; otherwise, $NX(b)$ and $SX(b)$ are arbitrarily
chosen different edges ending in $b$ when $X$ is $W$, and
beginning in $b$ when $X$ is $E$. (Note that for a non-transversal
inner vertex that belongs to an in-going tree there is no $E$-ward
branching, and for one that belongs to an out-going tree there is
no $W$-ward branching.)

It remains to define the values of the functions in $\L$ for the
transversal vertices, if there are such vertices in $D$. Let
$a_1,\ldots,a_n$, for $n\geq 2$, be our transversal of $D$. For
$i\in\{1,\ldots,n\mn 1\}$, if $(a_i,a_{i+1})$ is an edge of $D$,
then $SE(a_i)=NW(a_{i+1})=(a_i,a_{i+1})$, and if $(a_{i+1},a_i)$
is an edge of $D$, then $SW(a_i)=NE(a_{i+1})=(a_{i+1},a_i)$. (The
other possibility would be to take that if $(a_i,a_{i+1})$ is an
edge of $D$, then $NE(a_i)=SW(a_{i+1})=(a_i,a_{i+1})$, and if
$(a_{i+1},a_i)$ is an edge of $D$, then
$NW(a_i)=SE(a_{i+1})=(a_{i+1},a_i)$.)

For example, $a,b,c,d,e,f,g$ is the transversal of the K-graph
given below, and we define $SW(a)=NE(b)=(b,a)$,
$SW(b)=NE(c)=(c,b)$, $SE(c)=NW(d)=(c,d)$, etc., as it is suggested
by the following picture:
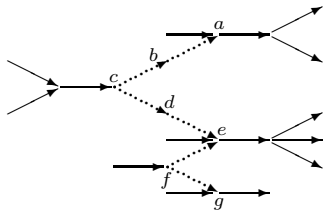
\begin{center}
\begin{picture}(130,90)
\put(0,40){\vector(2,1){19}} \put(0,60){\vector(2,-1){19}}
\put(20,50){\vector(1,0){19}} \put(40,20){\vector(1,0){19}}
\put(60,10){\vector(1,0){17.5}} \put(60,30){\vector(1,0){17.5}}
\put(60,70){\vector(1,0){17.5}} \put(80,10){\vector(1,0){19}}
\put(80,30){\vector(1,0){19}} \put(80,70){\vector(1,0){19}}
\put(100,30){\vector(1,0){19}} \put(100,29){\vector(2,-1){19}}
\put(100,31){\vector(2,1){19}} \put(100,69){\vector(2,-1){19}}
\put(100,71){\vector(2,1){19}}

\multiput(41,50)(2,1){18}{\makebox(0,0){\circle*{.5}}}
\put(74.5,67){\vector(2,1){5}} \put(55.5,57.5){\vector(2,1){5}}

\multiput(41,50)(2,-1){18}{\makebox(0,0){\circle*{.5}}}
\put(74.5,33){\vector(2,-1){5}} \put(55.5,42.5){\vector(2,-1){5}}

\multiput(61,20)(2,1){9}{\makebox(0,0){\circle*{.5}}}
\put(74.5,27){\vector(2,1){5}}

\multiput(61,20)(2,-1){9}{\makebox(0,0){\circle*{.5}}}
\put(74.5,13){\vector(2,-1){5}}

\put(80,8){\scriptsize\makebox(0,0)[t]{$g$}}
\put(60,18){\scriptsize\makebox(0,0)[t]{$f$}}
\put(81,35){\scriptsize\makebox(0,0)[t]{$e$}}
\put(61,47){\scriptsize\makebox(0,0)[t]{$d$}}
\put(40,55){\scriptsize\makebox(0,0)[t]{$c$}}
\put(55,65){\scriptsize\makebox(0,0)[t]{$b$}}
\put(80,75){\scriptsize\makebox(0,0)[t]{$a$}}

\end{picture}
\end{center}
The remaining values of the functions in $\L$ for transversal
vertices may be chosen freely provided we take care to guarantee
that $\langle D,\L\rangle$ separates $N$ from $S$.

It is clear that $\langle D,\L\rangle$ so defined separates $N$
from $S$. It remains to verify that every path of $D$ is decent in
$\langle D,\L\rangle$. If there were a path $b_1,\ldots,b_m$, with
$m\geq 2$, of $D$ that is not decent in $\langle D,\L\rangle$,
then all the vertices in this path would be transversal. This path
coincides either with $a_{j+1},\ldots,a_{j+m}$ or with
$a_{j+m},\ldots,a_{j+1}$, where $0\leq j$ and $j\! +\! m\leq n$. In
the first case, $(b_1,b_2)=SE(b_1)$, while
$(b_{m-1},b_m)=NW(b_m)$, which yields that the path
$b_1,\ldots,b_m$ is decent, contrary to our assumption. In the
second case, $(b_1,b_2)=NE(b_1)$, while $(b_{m-1},b_m)=SW(b_m)$,
which yields again that the path $b_1,\ldots,b_m$ is decent. So
every path of $D$ is decent in $\langle D,\L\rangle$. \qed

\section{Adding the identity graphs to K-graphs}
Our notion of K-graph, and the equivalent notions of global and
local K-graph, could be extended a little bit by allowing as
K-graphs oriented graphs of the form
\begin{center}
\begin{picture}(40,20)
\put(0,10){\vector(1,0){40}} \put(-3,10){\circle{2}}
\put(43,10){\circle{2}}
\end{picture}
\end{center}
with two vertices, one a $W$-vertex and the other an $E$-vertex;
these oriented graphs have a single edge made of these two
vertices, and they have no inner vertex. These additional K-graphs
would serve to represent identity deductions, which are related to
the sequents $A\vdash A$, and we will call them \emph{identity
graphs}.

For every oriented graph $D$ and every identity graph $I$ we will
have that the oriented graphs $D[-]I$ and $I[-]D$, with $[-]$
replaced by an appropriate $[e_W\mn e_E]$, are both equal to $D$ up to replacement of vertices.
The construction $I'$ of an identity graph $I$ would be a
single-node tree with $I$ in this unique node, and the
distinguished edges all being the unique edge of $I$. The
definition of construction involves now an appropriate
modification of (XYD). The notion of $\rho$-equivalence would be
extended so that for every construction $G$ we would have that
$G[-]I'$ is $\rho$-equivalent to $I'[-]G$, which is
$\rho$-equivalent to~$G$.

In the definition of local compass graph of Section~3 and in the
definition of K-graph of Section~5, in condition (3)
we would just replace the requirement that $D$ has an inner vertex by the requirement that it has an edge,
while everything else in these definitions would remain the same.

\section{Q-graphs}
If we determined the graphs produced by the rule (PC) of the
Introduction in the same manner as we determined in this paper the
graphs produced by planar plural cuts, we would obtain something
more general and more simple to characterize.

For the
definition of the new notion of global K-graph one possibility is
to reject in the definition of basic K-graph the requirement
(XYB). Everything else in the definition of construction and
global K-graph of Section~2 would remain unchanged. Let the new
global K-graphs be called \emph{global} Q-\emph{graphs}.

The new global Q-graphs can however be characterized more simply.
Let a Q-\emph{graph} be defined as an oriented graph $D$ that
satisfies conditions (1)-(3) of the definition of local compass
graph (see Section~3). The same three conditions are also found in
the definition of K-graph of Section~5. A notion of graph
associated with plural cuts in a context with the structural rule
of permutation, which, as our notion of Q-graph, is based
essentially on connectedness and non-circularity, may be found in
\cite{SS78}.

As for constructions of global K-graphs, we say that $G$ is a construction of a global Q-graph $D$ when $D$ is the root graph of $G$. One can show the following for every $X\in\{W,E\}$.

\prop{Proposition 7.1}{For every Q-graph $D$ and every $X$-edge
$d$ of $D$ there is a construction $G$ of $D$ such that $NX(G)=SX(G)=d$.}

\dkz We proceed by induction on the number $n$ of inner edges of
$D$. In the basis, when $n=0$, we rely on the new definition of
basic Q-graph (i.e.\ basic K-graph without (XYB)). In the induction
step we have $D=D_W[e_W\mn e_E]D_E$. By the induction hypothesis,
we have two constructions $G_X$ and $G_{\bar{X}}$ with root graphs
$D_X$ and $D_{\bar{X}}$ respectively such that if $d$ is in $D_X$,
then $NX(G_X)=SX(G_X)=d$ and
$NX(G_{\bar{X}})=SX(G_{\bar{X}})=e_{\bar{X}}$, and if $d$ is in
$D_{\bar{X}}$, then $NX(G_{\bar{X}})=SX(G_{\bar{X}})=d$ and
$N\bar{X}(G_X)=S\bar{X}(G_X)=e_X$. One can then verify that (XYC)
is satisfied, and that $NX(G)=SX(G)=d$, according to (XYD). \qed

As a corollary of this proposition we have that every Q-graph is a
global Q-graph. The converse being trivial, we have that the two
notions are equivalent.

This means that global Q-graphs could be defined by constructions
$G$ that do not involve at all the distinguished edges $YX(G)$.
For two arbitrary Q-graphs $D_W$ and $D_E$, an arbitrary $W$-edge
$e_W$ of $D_W$ and an arbitrary $E$-edge $e_E$ of $D_E$, the
oriented graph $D_W[e_W\mn e_E]D_E$ is a Q-graph. We need not
pay attention to (XYC) any more.

One could envisage the notion of Q-graph enlarged with identity
graphs, as in Section~6. The Q-graphs could be described in the
manner in which we have described K-graphs after Proposition 5.1,
which should still be applied (see also the Introduction). As a
K-graph, a Q-graph is made of a transversal and in-going and
out-going trees rooted in it. The difference is only that the
transversal need not be linear any more.

\vspace{3ex}

\noindent {\small \emph{Acknowledgement.} This work was supported by
the Ministry of Science of Serbia (Grant ON174026).}


\begin{thebibliography}{99}

\bibitem{A91} \textsc{V.M.\ Abrusci}, \textit{Phase semantics and sequent calculus for
pure noncommutative classical linear propositional logic},
\textbf{\textit{The Journal of Symbolic Logic}}, vol.\ 56 (1991),
pp.\ 1403-1451

\bibitem{B87} \textsc{B.\ B\' erard}, \textit{Formal properties of literal
shuffle}, \textbf{\textit{Acta Cybernetica}}, vol.~8 (1987), pp.\
27-39

\bibitem{CS97} \textsc{J.R.B.\ Cockett} and {\sc R.A.G.\ Seely}, \textit{Weakly
distributive categories}, \textbf{\textit{Journal of Pure and
Applied Algebra}}, vol.\ 114 (1997), pp.\ 133-173 (version with
some corrections at: http://www.math.mcgill.ca/rags)

\bibitem{D99} \textsc{K.\ Do\v sen}, \textit{On passing from singular to plural consequences}, \textbf{\textit{Logic at
Work: Essays Dedicated to the Memory of Helena Rasiowa}} (E.\
Or{\l}owska, editor), Physica-Verlag, Heidelberg, 1999, pp.\
533-547

\bibitem{G81} \textsc{D.M.\ Gabbay}, \textbf{\textit{Semantical Investigations in Heyting's Intuitionistic Logic}}, Reidel, Dordrecht, 1981

\bibitem{G35} \textsc{G.\ Gentzen}, \textit{Untersuchungen \" uber das
logische Schlie\ss en}, \textbf{\textit{Ma\-the\-mati\-sche
Zeit\-schrift}}, vol.\ 39 (1935), pp.\ 176-210, 405-431 (English
translation: \textit{Investigations into logical deduction},
\textbf{\textit{The Collected Papers of Gerhard Gentzen}}, M.E.\
Szabo, editor, North-Holland, Amsterdam, 1969, pp.\ 68-131)

\bibitem{H69} \textsc{F.\ Harary}, \textbf{\textit{Graph Theory}}, Addison-Wesley,
Reading, Mass., 1969

\bibitem{HS95} \textsc{J.\ Hudelmaier} and \textsc{P.\ Schroeder-Heister},
\textit{Classical Lambek logic}, \textbf{\textit{Theorem Proving with Analytic Tableaux
and Related Methods}} (P.\ Baumgartner et al., editors), Springer,
Berlin, 1995, pp.\ 247-262

\bibitem{K05} \textsc{J.\ Koslowski}, \textit{A monadic approach to polycategories},
\textbf{\textit{Theory and Applications of Categories}},
vol.\ 15 (2005), pp.\ 125-156

\bibitem{L89} \textsc{J.\ Lambek},
{\it Multicategories revisited}, \textbf{\textit{Categories in
Computer Science and Logic}} (J.W.\ Gray and A.\ Scedrov,
editors), American Mathematical Society, Providence, 1989, pp.\
217-239

\bibitem{L93} --------,
\textit{From categorial grammar to bilinear logic},
\textbf{\textit{Substructural Logics}} (K.\ Do\v sen and P.\
Schroeder-Heister, editors), Oxford University Press, Oxford, 1993,
pp.\ 207-237

\bibitem{S74} \textsc{D.S.\ Scott}, \textit{Completeness and axiomatizability in many-valued logic},
\textbf{\textit{Proceedings of the Tarski Symposium}} (L.\ Henkin
et al., editors), American Mathematical Society, Providence, 1974,
pp.\ 411-435

\bibitem{SS78} \textsc{D.J.\ Shoesmith} and {\sc T.J.\ Smiley},
\textbf{\textit{Multiple-Conclusion Logic}}, Cambridge University
Press, Cambridge, 1978

\bibitem{S75} \textsc{M.E.\ Szabo}, \textit{Polycategories},
\textbf{\textit{Communications in Algebra}}, vol.\ 3 (1975), pp.\
663-689

\end{thebibliography}
\end{document}